\documentclass[11pt]{article}
\usepackage[a4paper, margin=1in]{geometry}
\usepackage{amsfonts}
\usepackage{amsmath}
\usepackage{hyperref}
\usepackage{cite}
\usepackage{authblk}
\usepackage{titlesec}
\usepackage{pgfplots}
\pgfplotsset{compat=1.16}

\titlelabel{\thetitle.\quad}

\title{Initial Value Problem for a Caputo Space-time Fractional Schr{\"o}dinger Equation for the Delta Potential}
\author{Sepideh Saberhaghparvar\footnote{sepideh.saberhaghparvar65@gmail.com} {} and 
        Hossein Panahi\footnote{t-panahi@guilan.ac.ir (Corresponding author)}}
\affil{Department of Physics, Faculty of Science, University of Guilan, Rasht 51335-1914, Iran}

\begin{document}

\maketitle

\begin{abstract}
\noindent \textbf{Abstract:} In this paper, we investigate the initial value problem for a Caputo space-time fractional Schr{\" o}dinger equation for the delta potential. To solve this equation, we use the joint Laplace transform on the spatial coordinate and the Fourier transform on the time coordinate. Finally, the wave function and the time dependent energy eigenvalues are obtained for a particle which is subjected to the delta potential.

\vspace{5mm}

\noindent \textbf{keywords:} The fractional Schr{\"o}dinger equation; Caputo space-time fractional derivative

\vspace{5mm}

\noindent \textbf{PACS:} 03.65.Ca, 02.50.Ey, 02.30.Gp, 03.65.Db
\end{abstract}

\section{\label{section1}Introduction}
The fractional calculus is a generalization of the usual calculus, so that derivatives and integrals are defined for arbitrary real numbers. In some of the phenomena, the fractional operators simulate phenomena better than ordinary derivatives and normal integrals. The fractional calculus has been used in science and engineering. \cite{1,2,3} Recently the fractional Schr{\" o}dinger equation is studied in many fields, such as obstacle problem, phase transition and anomalous diffusion \cite{4,5,6,7,8,9,10} and etc. The fractional calculus began with Leibniz (1695-1697) and Euler's speculations (1730). After that, Riemann, Liouville, Gr{\" u}nwald and Letnikov \cite{11,12,13,14} provided definition of fractional derivatives. In 2000, Laskin started the first applications to quantum mechanics by considering the path integral formulation of quantum mechanics over L{\' e}vy paths and showed that the corresponding equation of motion is the space fractional Schr{\" o}dinger equation. \cite{15,16,17,18,19,20,21,22} The time fractional Schr{\" o}dinger equation discussed with Naber in 2004. \cite{23,24,25,26,27,28} He obtained the wave function of the time fractional Schr{\" o}dinger equation in terms of Mittag-Leffler function for a free particle and for a potential well. The time fractional Schr{\" o}dinger equation is obtained by replacing the time derivative of integer order with the derivative of the non-integer order. There are several ways for converting of the derivative of integer order to the derivative of the arbitrary real number for example Riemann-Liouville and Caputo definitions. In this paper, the Caputo fractional derivative has been used. The Caputo derivative of a function is defined as \cite{11}
\begin{eqnarray}
    {}_0^cD_t^\alpha f\left( t \right) = \frac{1}{\Gamma \left( n - \alpha  \right)} \int_0^t \frac{f^{\left( n \right)}\left( \tau \right) d\tau } {\left( {t - \tau } \right)^{\alpha  + 1 - n}}, \qquad \left( n - 1 \right) < \alpha  < n, \label{eq:1}
\end{eqnarray}
where for the case $\alpha \to n$, the Caputo derivative becomes an ordinary n-th derivative of the function $f\left( t \right)$.

In this paper, we consider the space-time fractional Schr{\" o}dinger equation for the delta potential, then we obtain the corresponding wave function by using the joint Laplace transform with respect to the time coordinate $t$ and the Fourier transform with respect to the spatial coordinate $x$ and by imposing one special initial condition. The Laplace transform of the Caputo derivative is given by \cite{11}
\begin{eqnarray}
    \mathcal{L}\left\{ {}_0^cD_t^\alpha f\left( t \right) \right\} = s^\alpha F\left( s \right) - \sum_{k = 0}^{n - 1} {s^{\alpha - k - 1} f^{\left( k \right)}\left( 0 \right)}, \label{eq:2}
\end{eqnarray}
where $F\left( s \right) = \int_0^\infty  {e^{-st}f\left( t \right) dt}$.

The paper is organized as follows. In section \ref{section2}, we introduce the space-time fractional Schr{\" o}dinger equation and one special initial condition. In section \ref{section3}, we obtain the wave function corresponding to the initial imposed condition in terms of Fox's H-function with two variables by applying the joint Laplace and Fourier transforms. In section \ref{section4}, we obtain the time dependent energy values by asymptotic expansion of Fox's H-function with two variables. The paper ends with the conclusions in section \ref{section5}.

\section{\label{section2}The space-time fractional Schr{\" o}dinger equation for delta potential}
The one-dimensional time dependent Schr{\" o}dinger equation is given by
\begin{eqnarray}
    i\hbar\frac{\partial \varphi \left( x,t \right)} {\partial t} = -\frac{\hbar^2}{2m}\frac{\partial ^2} {\partial x^2} \varphi \left( x,t \right) + V\left( x \right)\varphi \left( x,t \right), \label{eq:3}
\end{eqnarray}
where $m$ is the mass of the particle which is subjected to the potential $V\left( x \right)$ and $\hbar$ is the Plank's constant. Also the one-dimensional space-time fractional Schr{\" o}dinger equation is written as follows \cite{23}
\begin{eqnarray}
    i^\alpha \eta {}_0^cD_t^\alpha \varphi \left( x,t \right) = - \frac{1}{2}mc^2 \left( \frac{\hbar}{mc} \right)^{2\beta} {}_0^cD_x^{2\beta}\varphi \left( x,t \right) + V\left( x \right)\varphi \left( x,t \right), \label{eq:4}
\end{eqnarray}
where $0 < \alpha \le 1$, $0.5 < \beta \le 1$, $\eta = mc^2 \left( \frac{\hbar}{mc^2} \right)^\alpha$ and the scaling factors $i^\alpha \eta$ and $-\frac{1}{2}mc^2 \left( \frac{\hbar}{mc} \right)^{2\beta}$ have been added to equalize units on both sides of the equation (\ref{eq:4}). It is easy to see that for the special case $\alpha=1$ and $\beta=1$, the space-time fractional Schr{\" o}dinger equation (\ref{eq:4}) reduces to the ordinary Schr{\" o}dinger equation (\ref{eq:3}). For delta potential $V\left( x \right) = -v_0 \delta \left( x \right)$, we have
\begin{eqnarray}
    i^\alpha \eta {}_0^cD_t^\alpha \varphi \left( x,t \right) =  - \frac{1}{2} mc^2 \left( \frac{\hbar}{mc} \right)^{2\beta } {}_0^cD_x^{2\beta} \varphi \left( x,t \right) - v_0 \delta \left( x \right)\varphi \left( x,t \right). \label{eq:5}
\end{eqnarray}
We now consider two physical boundary and initial conditions as
\begin{eqnarray}
    \begin{cases}
    \varphi \left( x,0 \right) = g\left( x \right), \\
    \varphi \left( x,t \right) \to 0, & as\quad\left| x \right| \to \infty. \\
    \varphi \left( 0,t \right) = 0, & n \ge 1 \\
    \varphi _x^{\left( n \right)}\left( x,0 \right) = 0,
    \end{cases} \label{eq:6}
\end{eqnarray}
To solve the Eq. (\ref{eq:5}), we apply the joint Laplace transform on the spatial coordinate and the Fourier transform on the time coordinate defined by \cite{29}
\begin{eqnarray}
    \overline {\tilde \varphi } \left( k,s \right) = \frac{1}{\sqrt {2\pi }} \int_{ - \infty }^{ + \infty } {e^{-ikx} dx} \int_0^\infty  {e^{-st} \varphi \left( x,t \right) dt}, \label{eq:7}
\end{eqnarray}
where the signs of ($-$) and ($\sim$) are used to denote the Laplace and the Fourier transforms respectively, also $k$ and $s$ are the Fourier and the Laplace transform variables respectively.

\section[Wave function of one-dimensional space-time fractional Schr{\" o}dinger equation for the initial imposed condition]{\label{section3}Wave function of one-dimensional space-time fractional \\ Schr{\" o}dinger equation for the initial imposed condition}
By applying the joint Laplace and Fourier transforms on Eq. (\ref{eq:5}) and using the boundary conditions of Eq. (\ref{eq:6}), we obtain the following equation
\begin{eqnarray}
    \overline {\tilde \varphi } \left( k,s \right) = s^{\alpha - 1} \frac{\tilde g \left( k \right)} {s^\alpha + \frac{F \left( ik \right)^{2\beta}} {i^\alpha \eta }}, \label{eq:8}
\end{eqnarray}
where $F = \frac{1}{2} mc^2 \left( \frac{\hbar}{mc} \right)^{2\beta }$. The inverse Laplace transform of Eq. (\ref{eq:8}) gives
\begin{eqnarray}
    \tilde \varphi \left( k,t \right) = \tilde g \left( k \right) E_\alpha \left( - \frac{F \left( ik \right)^{2\beta } t^\alpha } {i^\alpha \eta } \right). \label{eq:9}
\end{eqnarray}
In the above calculation, we have used the following formula \cite{30}
\begin{eqnarray}
    \mathcal{L}^{-1} \left\{\frac{m! s^{\alpha  - \beta }} {\left( s^\alpha \pm a^2 \right)^{m + 1}} \right\} = t^{\alpha m + \beta - 1} E_{\alpha ,\beta }^{\left( m \right)}\left( \mp a^2{t^\alpha } \right), \label{eq:10}
\end{eqnarray}
where $E_{\alpha ,\beta }\left( z \right)$ is the Mittag-Leffler function defined as \cite{31}
\begin{eqnarray}
    E_{\alpha ,\beta } \left( z \right) = \sum_{n = 0}^\infty \frac{z^n}{\Gamma\left( n \alpha  + \beta \right)}, \qquad \alpha > 0, \beta \in \mathbb{C}. \label{eq:11}
\end{eqnarray}
Now, we assume the function $g\left(x\right)$ as
\begin{eqnarray}
    g\left( x \right) = A E_{\alpha '} \left(  - \lambda \left| x \right|^{\alpha '} \right), \qquad \alpha ' \in \mathbb{C}. \label{eq:12}
\end{eqnarray}
Hence in Eq. (\ref{eq:9}), we must use its Fourier transform as $\tilde g \left( k \right)$, so by applying the Fourier transform to the relation (\ref{eq:12}), we have
\begin{eqnarray}
    \tilde g \left( k \right)
    && = \frac{A}{\sqrt{2\pi}} \int_{ - \infty }^{ + \infty } { e^{ - ikx} E_{\alpha '} \left( - \lambda {\left| x \right|}^{\alpha '} \right) dx } \nonumber \\
    && = \frac{2A}{\sqrt{2\pi}} \int_0^\infty { \cos\left( kx \right) H_{1,2}^{1,1} \left[ \lambda x^{\alpha '} \left| 
    \begin{array}{c}
        \left( 0,1 \right) \\
        \left( 0,1 \right), \left( 0,\alpha' \right)
    \end{array} \right. \right] dx }, \label{eq:13}
\end{eqnarray}
where $H_{p,q}^{m,n}\left[ {x\left| - \right.} \right]$ denotes Fox's H-function as follows \cite{32}
\begin{eqnarray}
    H\left[ x \right] && = H_{p,q}^{m,n}\left[ {x\left|
    \begin{array}{c}
    {\left( a_j,\alpha_j \right)}_{1,p} \\
    {\left( b_j,\beta_j \right)}_{1,q}
    \end{array} \right.} \right] = H_{p,q}^{m,n}\left[ {x\left| 
    \begin{array}{c}
    \left( a_1,\alpha_1 \right), \ldots \left( a_p,\alpha_p \right) \\
    \left( b_1,\beta_1 \right), \ldots ,\left( b_q,\beta_q \right)
    \end{array} \right.} \right] \nonumber \\
    && = \frac{1}{2\pi i} \int_{\mathcal{L}} {\Theta\left( s \right) x^s ds} \nonumber \\ 
    && = \frac{1}{2\pi i} \int_{\mathcal{L}} {\frac{\prod \nolimits_{j = 1}^m \Gamma \left( b_j - \beta_j s \right) \prod \nolimits_{j = 1}^n \Gamma \left( 1 - a_j + \alpha_j s \right)} { \prod \nolimits_{j = m + 1}^q \Gamma \left( 1 - b_j + \beta _j s \right) \prod \nolimits_{j = n + 1}^p \Gamma \left( a_j - \alpha_j s \right)} x^s ds}. \label{eq:14}
\end{eqnarray}
In the above relation, we have used the following formula \cite{32}
\begin{eqnarray}
    H_{1,2}^{1,1}\left[ -z \left| 
    \begin{array}{c}
        \left( {0,1} \right) \\
        \left( {0,1} \right), \left( {0,\alpha'} \right)
    \end{array} 
    \right. \right] = E_{\alpha'} \left( z \right). \label{eq:15}
\end{eqnarray}
Now we use the following equation to solve the integral in Eq. (\ref{eq:13}) \cite{32}
\begin{eqnarray}
    && \int_0^\infty {x^{\rho - 1} \cos \left( ax \right) H_{p,q}^{m,n}\left[ {b x^\sigma \left| \begin{array}{c}
    \left( a_p,A_p \right) \\
    \left( b_q,B_q \right)
    \end{array} \right.} \right] dx} \nonumber \\
    && = \frac{2^{\rho  - 1} \sqrt \pi} {a^\rho } H_{p + 2,q}^{m,n + 1}\left[ {b{{\left( {\frac{2}{a}} \right)}^\sigma }\left| {\begin{array}{c}
    \left( {\frac{2 - \rho}{2},\frac{\sigma }{2}} \right),\left( {a_p,A_p} \right),\left( {\frac{{1 - \rho }}{2},\frac{\sigma }{2}} \right) \\
    \left( {{b_q},{B_q}} \right)
    \end{array}} \right.} \right], \label{eq:16}
\end{eqnarray}
where $a,\Delta ,\sigma > 0$, $\rho ,b \in \mathbb{C}$; $\left| {\arg b} \right| < \frac{1}{2}\pi \Delta$;
\begin{eqnarray}
    \Re \left( \rho \right) + \sigma \min_{1 \le j \le m} {\Re \left( \frac{b_j}{B_j} \right) } > 0; \quad \Re \left( \rho \right) + \sigma \max_{1 \le j \le n} \left[ \frac{\left( {a_j - 1} \right)}{A_j} \right] < 1, \label{eq:17}
\end{eqnarray}
and
\begin{eqnarray}
    \Delta  = \sum_{j = 1}^n {A_j} - \sum_{j = n + 1}^p {A_j} + \sum_{j = 1}^m {B_j} - \sum_{j = m + 1}^q {B_j}. \label{eq:18}
\end{eqnarray}
Thus the integral conditions for Eq. (\ref{eq:13}) are as
\begin{eqnarray}
    \alpha > 0, \alpha' < 2, \arg{\lambda} < \frac{\pi}{2} \left( 2 - \alpha' \right). \label{eq:19}
\end{eqnarray}
By using of (\ref{eq:16}) and applying in Eq. (\ref{eq:13}), we get
\begin{eqnarray}
    \tilde g \left( k \right) = \frac{\sqrt 2 A}{\left| k \right|} H_{3,2}^{1,2}\left[ \lambda {\left( \frac{2}{\left| k \right|} \right)}^{\alpha'} \left| 
    \begin{array}{c}
    \left( \frac{1}{2}, \frac{\alpha'}{2} \right),\left( 0,1 \right),\left( 0,\frac{\alpha'}{2} \right) \\
    \left( 0,1 \right),\left( {0,\alpha'} \right)
    \end{array} \right. \right]. \label{eq:20}
\end{eqnarray}
We can also use the properties of Fox's H-function given in Ref. \cite{32} and so we have
\begin{eqnarray}
    \tilde g \left( k \right) = \frac{\sqrt 2 A}{2\lambda^{\frac{1}{\alpha'}}} H_{2,3}^{2,1}\left[ {\frac{1}{\lambda } \left( \frac{\left| k \right|}{2} \right)^{\alpha '}\left| \begin{array}{c}
    \left( 1 - \frac{1}{\alpha'},1 \right),\left( 0,\alpha' \right) \\
    \left( 0,\frac{\alpha'}{2} \right),\left( {1 - \frac{1}{\alpha'},1} \right),\left( \frac{1}{2},\frac{\alpha '}{2} \right)
    \end{array} \right.} \right]. \label{eq:21}
\end{eqnarray}
Now by substituting (\ref{eq:21}) into (\ref{eq:9}), we get
\begin{eqnarray}
    \tilde \varphi \left( k,t \right) = && \frac{\sqrt 2 A} {2\lambda^{\frac{1}{\alpha'}}} H_{2,3}^{2,1} \left[ \frac{1}{\lambda}\left( \frac{\left| k \right|}{2} \right)^{\alpha'} \left| 
    \begin{array}{c}
    \left( 1 - \frac{1}{\alpha'},1 \right), \left( 0,\alpha' \right) \\
    \left( 0,\frac{\alpha'}{2} \right), \left( 1 - \frac{1}{\alpha'},1 \right), \left( \frac{1}{2},\frac{\alpha'}{2} \right)
    \end{array} \right. \right] \nonumber \\
    && \times E_\alpha \left( - \frac{F\left( ik \right)^{2\beta} t^\alpha} {i^\alpha \eta } \right), \label{eq:22}
\end{eqnarray}
where by calculating the inverse Fourier transform of (\ref{eq:22}) and using the integral containing of two Fox's H-function together applying the properties of Fox's H-function \cite{32,33}, we get the wave function $\varphi(x,t)$ as
\begin{eqnarray}
    && \varphi \left( x,t \right) 
    = \frac{2A\left| x \right|^{-1}} {\alpha'\beta \lambda^{\frac{1}{\alpha'}}} \nonumber \\
    && \times H_{2,0:1,2;2,3}^{0,1:1,1;2,1}\left[ \left. 
    \begin{array}{c}
    \frac{ -4 \left( \frac{F}{i^\alpha \eta } \right)^{\frac{1}{\beta }} t^{\frac{\alpha}{\beta }}} {\left| x \right|^2} \\
    \frac{\lambda^{-\frac{2}{\alpha'}}} {\left| x \right|^2}
    \end{array} \right|
    \begin{array}{c}
    \left( \frac{1}{2};1,1 \right),\left( 0;1,1 \right):\left( 0,\frac{1}{\beta} \right);\left( 1 - \frac{1}{\alpha'},\frac{2}{\alpha'} \right),\left( 0,2 \right) \\
    - :\left( 0,\frac{1}{\beta} \right),\left( 0,\frac{\alpha}{\beta} \right);\left( 0,1 \right),\left( 1 - \frac{1}{\alpha'},\frac{2}{\alpha'} \right),\left( \frac{1}{2},1 \right)
    \end{array} \right], \label{eq:23}
\end{eqnarray}
where $\beta > - \frac{1}{2}$, $\frac{\alpha}{\beta} \ge 0$, $2-\alpha > \beta$, $\alpha' < \frac{4}{3}$ and 
\begin{eqnarray*}
    H_{p_1,q_1:p_2,q_2;p_3,q_3}^{0,n_1:m_2,n_2;m_3,n_3} \left[ {\left. 
    \begin{array}{c} x \\ y \end{array}
    \right| 
    \begin{array}{c}
    {{{\left( {{a_i};{\alpha _i},{A_i}} \right)}_{1,{p_1}}}:{{\left( {{c_i},{\gamma_i}} \right)}_{1,{p_2}}};{{\left( {{e_i},{E_i}} \right)}_{1,{p_3}}}} \\
    {{{\left( {{b_j};{\beta _j},{B_j}} \right)}_{1,{q_1}}}:{{\left( {{d_j},{\delta_j}} \right)}_{1,{q_2}}};{{\left( {{f_j},{F_j}} \right)}_{1,{q_3}}}}
    \end{array}} \right]
\end{eqnarray*}
is the Fox's H-function of two variables which is defined as \cite{34}
\begin{eqnarray}
    H\left[ x,y \right] 
    && = H\left[ \begin{array}{c} x \\ y \end{array} \right] \nonumber \\
    && = H_{p_1,q_1:p_2,q_2;p_3,q_3}^{0,n_1:m_2,n_2;m_3,n_3}\left[ \left. 
    \begin{array}{c} x \\ y \end{array} \right|
    \begin{array}{c}
    \left( a_i;\alpha_i,A_i \right)_{1,p_1}:\left( c_i,\gamma_i \right)_{1,p_2};\left( e_i,E_i \right)_{1,p_3} \\
    \left( b_j;\beta_j,B_j \right)_{1,q_1}:\left( d_j,\delta_j \right)_{1,q_2};\left( f_j,F_j \right)_{1,q_3}
    \end{array} \right] \nonumber \\
    && =  - \frac{1}{4 \pi^2} \int_{\mathcal{L}_1}{ \int_{\mathcal{L}_2} {\phi _1 \left(\ \xi ,\eta \right) \theta _2 \left( \xi \right) \theta _3 \left( \eta \right) x^\xi y^\eta d\xi} d\eta}, \label{eq:24}
\end{eqnarray}
where
\begin{eqnarray}
    \phi_1\left( \xi ,\eta \right) && = \frac{ \prod_{j = 1}^{n_1} \Gamma \left( 1 - a_j + \alpha_j \xi + A_j \eta \right)} { \prod_{j = n_1 + 1}^{p_1} \Gamma \left( a_j - \alpha_j \xi - A_j \eta \right) \prod_{j = 1}^{q_1} \Gamma \left( 1 - b_j + \beta _j \xi + B_j\eta \right)}, \label{eq:25} \\
    \theta_2 \left( \xi  \right) && = \frac{ \prod_{j = 1}^{n_2} \Gamma \left( 1 - c_j + \gamma_j \xi \right) \prod_{j = 1}^{m_2} \Gamma \left( d_j - \delta _j \xi \right)} { \prod_{j = n_2 + 1}^{p_2} \Gamma \left( c_j - \gamma_j \xi \right) \prod_{j = m_2 + 1}^{q_2} \Gamma \left( 1 - d_j + \delta_j \xi \right)}, \label{eq:26} \\
    \theta_3 \left( \eta  \right) && = \frac{ \prod_{j = 1}^{n_3} \Gamma \left( 1 - e_j + E_j \eta  \right) \prod_{j = 1}^{m_3} \Gamma \left( f_j - F_j\eta  \right)} {\prod_{j = n_3 + 1}^{p_3} \Gamma \left( e_j - E_j \eta \right) \prod_{j = m_3 + 1}^{q_3} \Gamma \left( 1 - f_j + F_j \eta \right)}. \label{eq:27}
\end{eqnarray}
Before calculating the energy eigenvalues, we obtain the wave function $\varphi(x,t)$ in limits of $\alpha \to 1$, $\beta \to 1$ and $\alpha' \to 1$ cases.
\begin{eqnarray}
    \varphi \left( x,t \right) _{ \left| \begin{array}{c}
    \alpha \to 1 \\
    \beta \to 1 \\
    \alpha' \to 1
    \end{array} \right.} = \frac{2A\left| x \right|^{-1}} {\lambda } H_{2,0:0,1;1,1}^{0,1:1,1;1,1}\left[ \left. \begin{array}{c}
    \frac{-4Ft}{i \eta \left| x \right|^2} \\
    \frac{\lambda^{-2}} {\left| x \right|^2}
    \end{array} \right|\begin{array}{c}
    \left( \frac{1}{2};1,1 \right),\left( 0;1,1 \right): - ;\left( 0,1 \right) \\
    { - :\left( 0,1 \right);\left( 0,1 \right)}
    \end{array} \right], \label{eq:28}
\end{eqnarray}
where by using of Eq.(\ref{eq:24}) we get
\begin{eqnarray}
    \varphi \left( x,t \right) _{\left| \begin{array}{c}
    \alpha \to 1 \\
    \beta \to 1 \\
    \alpha' \to 1
    \end{array} \right.} 
    = && \frac{2A\left| x \right|^{-1}} {\lambda \sqrt \pi} \left( \frac{1}{2\pi i} \right) \sum_{n = 0}^\infty \frac{\left( {-1} \right)^n} {n!} {\left( \frac{-4Ft} {i\eta \left| x \right|^2} \right)^n} \nonumber \\
    && \times \int_{L_1} \frac{\Gamma \left( -\eta' \right) \Gamma \left( 1 + \eta' \right) \Gamma \left( \frac{1}{2} + n + \eta' \right)} {\Gamma \left( { - n - \eta'} \right)} \left( 4\frac{\lambda^{-2}} {\left| x \right|^2} \right)^{\eta'}d\eta'. \label{eq:29}
\end{eqnarray}
Using the properties of Gamma functions as
\begin{eqnarray}
    \Gamma \left( \frac{1}{2} + n + \eta' \right) = \left( {\frac{1}{2} + \eta '} \right)_n \Gamma \left( \frac{1}{2} + \eta' \right), \label{eq:30} \\
    \Gamma \left( -n-\eta' \right) = \left( -\eta' \right)_{-n} \Gamma \left( -\eta' \right) = \frac{\left( -1 \right)^n} {\left( 1+\eta' \right)_n} \Gamma \left( -\eta' \right), \label{eq:31}
\end{eqnarray}
and substituting them into Eq. (\ref{eq:29}), we have
\begin{eqnarray}
    \varphi \left( x,t \right) _{\left| \begin{array}{c}
    \alpha  \to 1 \\
    \beta  \to 1 \\
    \alpha' \to 1
    \end{array} \right.} = && \frac{2A\left| x \right|^{-1}} {\lambda \sqrt \pi  } \sum_{n = 0}^\infty  \frac{1}{n!} \left( \frac{-4Ft} {i\eta \left| x \right|^2} \right)^n \nonumber \\
    && \times H_{2,0}^{0,2}\left[ \left. \frac{4}{\left| x \right|^2{\lambda^2}} \right| \begin{array}{c}
    \left( \frac{1}{2}-n,1 \right),\left( -n,1 \right) \\
    - 
    \end{array} \right]. \label{eq:32}
\end{eqnarray}
If we continue the computational process of Fox's H-function according to Ref. \cite{34}, we get
\begin{eqnarray}
    \varphi \left( x,t \right) _{\left| \begin{array}{c}
    \alpha \to 1 \\
    \beta \to 1 \\
    \alpha' \to 1
    \end{array} \right.} = A e^{\frac{iFt{\lambda^2}}{\eta}} e^{ -\lambda \left| x \right|} = A e^{\frac{i t\lambda^2\hbar}{2m}} e^{ -\lambda \left| x \right|}, \label{eq:33}
\end{eqnarray}
where up to a coefficient A, it is the wave function of ordinary Schr{\" o}dinger equation for delta potential
\begin{eqnarray}
    \varphi \left( x,t \right) = \frac{\sqrt {v_0 m} } {\hbar} e^{ -\frac{\left| x \right|v_0 m} {\hbar^2}} e^{i\frac{v_0^2m}{2\hbar^3}t}. \label{eq:34}
\end{eqnarray}
Now by comparing (\ref{eq:33}) with (\ref{eq:34}) and regardless of the normalization coefficient, the parameter must be as $\lambda = \frac{v_0 m}{\hbar^2}$ and $A = \frac{\sqrt {v_0 m}}{\hbar}$.

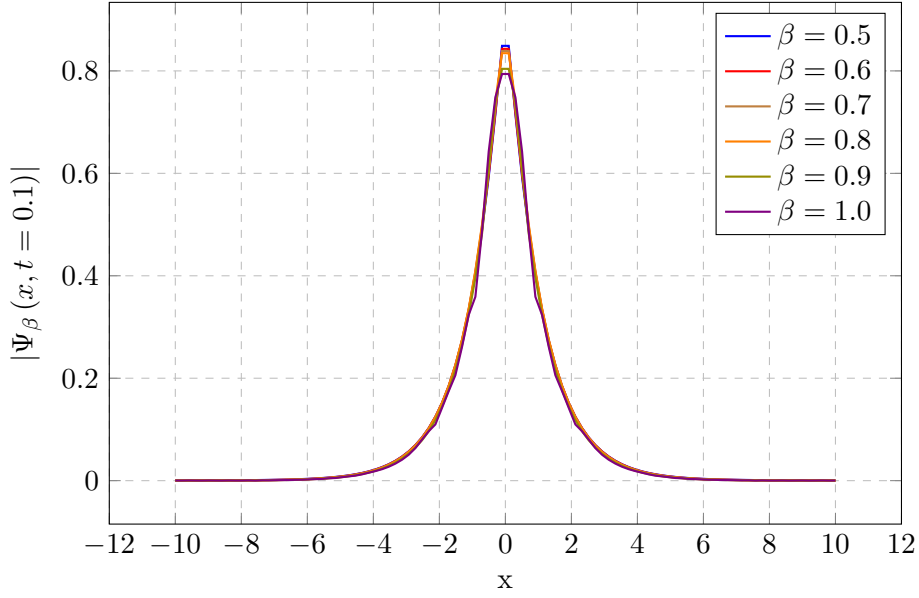
\begin{figure}
\centering
    \begin{tikzpicture}
        \begin{axis}[no markers,grid=major,width=12cm,height=8.5cm,grid style=dashed,
            xlabel = {x},
            ylabel = {{$\left|\Psi_\beta\left(x,t=0.1\right)\right|$}}]
            \addplot[thick,solid,blue] table [x = x, y = y1, col sep = comma] {plot1.csv};
            \addlegendentry{$\beta=0.5$}
            \addplot[thick,solid,red] table [x = x, y = y2, col sep = comma] {plot1.csv};
            \addlegendentry{$\beta=0.6$}
            \addplot[thick,solid,brown] table [x = x, y = y3, col sep = comma] {plot1.csv};
            \addlegendentry{$\beta=0.7$}
            \addplot[thick,solid,orange] table [x = x, y = y4, col sep = comma] {plot1.csv};
            \addlegendentry{$\beta=0.8$}
            \addplot[thick,solid,olive] table [x = x, y = y5, col sep = comma] {plot1.csv};
            \addlegendentry{$\beta=0.9$}
            \addplot[thick,solid,violet] table [x = x, y = y6, col sep = comma] {plot1.csv};
            \addlegendentry{$\beta=1.0$}
        \end{axis}
    \end{tikzpicture}
    \caption{\label{fig:1}Plot of Eq. (\ref{eq:23}) in terms of $x$ for $\alpha=1$, $\alpha'=1$, $t=0.1$ and different values of $\beta$ and for simplicity $m=\hbar=c=\lambda=1$.}
\end{figure}
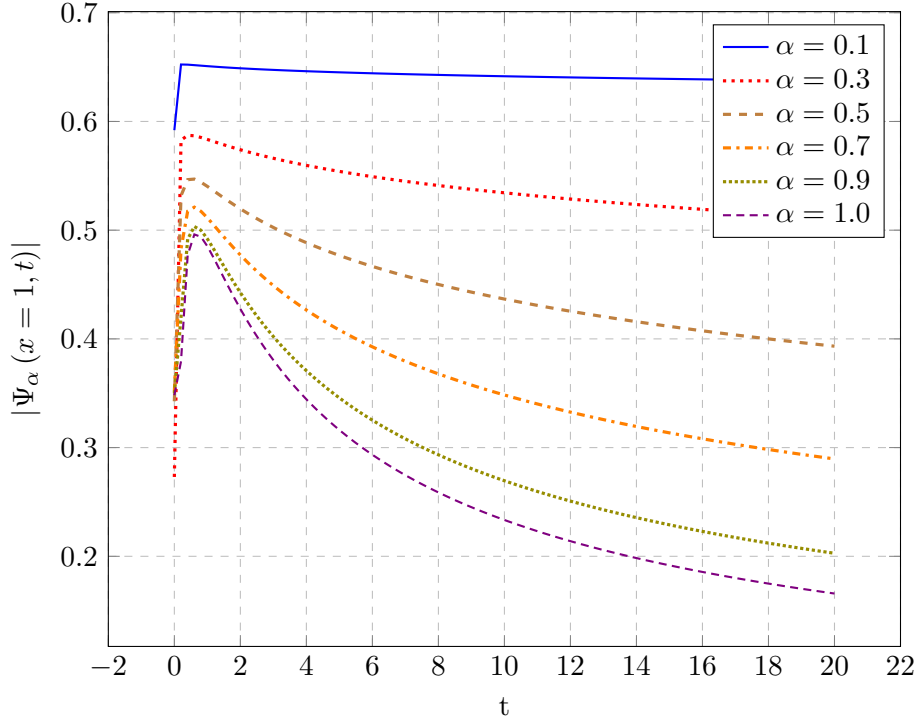
\begin{figure}
\centering
    \begin{tikzpicture}
        \begin{axis}[width = 12cm, height = 10cm, no markers, grid=major, grid style=dashed,
            xlabel = {t},
            ylabel = {{$\left|\Psi_\alpha\left(x=1,t\right)\right|$}}]
            \addplot[thick,solid,blue] table [x = t, y = y1, col sep = comma] {plot2.csv};
            \addlegendentry{$\alpha=0.1$}
            \addplot[very thick,dotted,red] table [x = t, y = y2, col sep = comma] {plot2.csv};
            \addlegendentry{$\alpha=0.3$}
            \addplot[very thick,dashed,brown] table [x = t, y = y3, col sep = comma] {plot2.csv};
            \addlegendentry{$\alpha=0.5$}
            \addplot[very thick,dashdotted,orange] table [x = t, y = y4, col sep = comma] {plot2.csv};
            \addlegendentry{$\alpha=0.7$}
            \addplot[very thick,densely dotted,olive] table [x = t, y = y5, col sep = comma] {plot2.csv};
            \addlegendentry{$\alpha=0.9$}
            \addplot[thick,densely dashed,violet] table [x = t, y = y6, col sep = comma] {plot2.csv};
            \addlegendentry{$\alpha=1.0$}
        \end{axis}
    \end{tikzpicture}    
    \caption{\label{fig:2}Plot of Eq. (\ref{eq:23}) in terms of $t$ for $\beta=1$, $\alpha'=1$, $x=1$ and some values of $\alpha$ and for simplicity $m=\hbar=c=\lambda=1$.}
\end{figure} 

In Fig. \ref{fig:1} we have plotted the space fractional wave function (\ref{eq:23}) for some values of $\beta$ by considering $\alpha=1$, $\alpha'=1$ and Fig. \ref{fig:2} shows the time fractional wave function for different values of $\alpha$ by considering $\beta=1$, $\alpha'=1$. In both figures the bivariate $H$-function was computed using a Python code provided in Ref. \cite{35} and available at https://github.com/melayadi/multivarFoxH.

According to Fig. \ref{fig:1}, the space wave function is consistent with the wave function of the ordinary Schr{\" o}dinger equation for $\beta=1$ for a given time ($t=0.1$). Also, Fig. \ref{fig:2} shows that the wave function decreases exponentially with increasing time in the case of $\alpha=1$ for a specific location ($x=1$) which is consistent with the wave function of the ordinary Schr{\" o}dinger equation. Thus, we have obtained the time fractional Schr{\" o}dinger equation for the delta potential by assuming the initial function (\ref{eq:12}). In the next section, we obtained the energy eigenvalues of the system.

\section{\label{section4}The energy eigenvalues of delta potential for the initial imposed condition}
For the space-time fractional Schr{\" o}dinger equation, the energy eigenvalues are as \cite{35}
\begin{eqnarray}
    E_\alpha  = i^\alpha \eta \int_{ -\infty }^{ +\infty } \varphi^* \left( x,t \right) {}_0^cD_t^\alpha \varphi \left( x,t \right) dx. \label{eq:35}
\end{eqnarray}
The wave function $\varphi(x,t)$ in Eq. (\ref{eq:23}) can be rewritten using Fox's H-function identities as follows
\begin{eqnarray}
    && \varphi \left( {x,t} \right) 
    = \frac{Ai} {\alpha' \beta \left( \lambda \right)^{\frac{1}{\alpha'}} \left( \frac{F}{i^\alpha \eta } \right)^{\frac{1}{2\beta }} t^{\frac{\alpha}{2\beta}}} \nonumber \\
    && \times H_{2,0:1,2;2,3}^{0,1:1,1;2,1}\left[ \left. 
    \begin{array}{c}
    \frac{-4\left( \frac{F}{i^\alpha \eta } \right)^{\frac{1}{\beta}} t^{\frac{\alpha}{\beta}}} {\left| x \right|^2} \\
    \frac{\lambda^{-\frac{2}{\alpha'}}} {\left| x \right|^2}
    \end{array} \right|
    \begin{array}{c}
    \left( 1;1,1 \right), \left( \frac{1}{2};1,1 \right): \left( \frac{1}{2\beta}, \frac{1}{\beta } \right); \left( {1 - \frac{1}{{\alpha '}},\frac{2}{{\alpha '}}} \right),\left( {0,2} \right) \\
    - :\left( \frac{1}{2\beta}, \frac{1}{\beta} \right), \left( \frac{\alpha}{2\beta}, \frac{\alpha}{\beta} \right); \left( 0,1 \right), \left( 1-\frac{1}{\alpha'}, \frac{2}{\alpha'} \right),\left( \frac{1}{2},1 \right)
    \end{array} \right]. \label{eq:36}
\end{eqnarray}
On the other hand, for calculating of the Caputo fractional derivative of Fox's H-function of two variables from Eq. (\ref{eq:1}) for $n \to 1$, we must obtain the derivative of the first order of  Fox's H-function of two variables and substitute it in Eq. (\ref{eq:1}). Thus, the derivative of first order Fox's H-function of two variables is as \cite{34}
\small
\begin{eqnarray}
    && \frac{d\varphi \left( x,t \right)} {dt}
    = \frac{Ai} {\alpha' \beta \left( \lambda \right)^{\frac{1}{\alpha'}} \left( \frac{F}{i^\alpha \eta } \right)^{\frac{1}{2\beta }} t^{\frac{\alpha}{2\beta}+1}} \nonumber \\
    && \times H_{2,0:1,2;2,3}^{0,1:1,1;2,1}\left[ \left. 
    \begin{array}{c}
    \frac{-4\left( \frac{F}{i^\alpha \eta } \right)^{\frac{1}{\beta}} t^{\frac{\alpha}{\beta}}} {\left| x \right|^2} \\
    \frac{\lambda^{-\frac{2}{\alpha'}}} {\left| x \right|^2}
    \end{array} \right|
    \begin{array}{c}
    \left( 1;1,1 \right), \left( \frac{1}{2};1,1 \right): \left( \frac{1}{2\beta}, \frac{1}{\beta } \right); \left( {1 - \frac{1}{{\alpha '}},\frac{2}{{\alpha '}}} \right),\left( {0,2} \right) \\
    - :\left( \frac{1}{2\beta}, \frac{1}{\beta} \right), \left( 1+\frac{\alpha}{2\beta}, \frac{\alpha}{\beta} \right); \left( 0,1 \right), \left( 1-\frac{1}{\alpha'}, \frac{2}{\alpha'} \right),\left( \frac{1}{2},1 \right)
    \end{array} \right], \label{eq:37}
\end{eqnarray}
\normalsize
and so, we have
\small
\begin{eqnarray}
    && {}_0^cD_t^\alpha \varphi \left( x,t \right) 
    = \frac{1}{\Gamma \left( 1-\alpha \right)} \int_0^t \frac{Ai} {\alpha' \beta \left( \lambda \right)^{\frac{1}{\alpha'}} \left( \frac{F}{i^\alpha \eta } \right)^{\frac{1}{2\beta }} \tau^{\frac{\alpha}{2\beta}+1}} \nonumber \\
    && \times H_{2,0:1,2;2,3}^{0,1:1,1;2,1}\left[ \left. 
    \begin{array}{c}
    \frac{-4\left( \frac{F}{i^\alpha \eta } \right)^{\frac{1}{\beta}} \tau^{\frac{\alpha}{\beta}}} {\left| x \right|^2} \\
    \frac{\lambda^{-\frac{2}{\alpha'}}} {\left| x \right|^2}
    \end{array} \right| \right. \left.
    \begin{array}{c}
    \left( 1;1,1 \right), \left( \frac{1}{2};1,1 \right): \left( \frac{1}{2\beta}, \frac{1}{\beta} \right); \left( {1 - \frac{1}{\alpha'},\frac{2}{\alpha'}} \right),\left( 0,2 \right) \\
    - :\left( \frac{1}{2\beta}, \frac{1}{\beta} \right), \left( 1+\frac{\alpha}{2\beta}, \frac{\alpha}{\beta} \right); \left( 0,1 \right), \left( 1-\frac{1}{\alpha'}, \frac{2}{\alpha'} \right),\left( \frac{1}{2},1 \right)
    \end{array} \right] \nonumber \\
    && \times \left( t-\tau \right)^{-\alpha} d\tau = \frac{Ai} {\alpha' \beta \left( \lambda \right)^{\frac{1}{\alpha'}} \left( \frac{F}{i^\alpha \eta } \right)^{\frac{1}{2\beta }} t^{\frac{\alpha}{2\beta}+1}} \nonumber \\
    && \times H_{2,0:1,2;2,3}^{0,1:1,1;2,1}\left[ \left. 
    \begin{array}{c}
    \frac{-4\left( \frac{F}{i^\alpha \eta } \right)^{\frac{1}{\beta}} t^{\frac{\alpha}{\beta}}} {\left| x \right|^2} \\
    \frac{\lambda^{-\frac{2}{\alpha'}}} {\left| x \right|^2}
    \end{array} \right| \right. \left.
    \begin{array}{c}
    \left( 1;1,1 \right), \left( \frac{1}{2};1,1 \right): \left( \frac{1}{2\beta}, \frac{1}{\beta} \right); \left( {1 - \frac{1}{\alpha'},\frac{2}{\alpha'}} \right),\left( 0,2 \right) \\
    - :\left( \frac{1}{2\beta}, \frac{1}{\beta} \right), \left( \alpha+\frac{\alpha}{2\beta}, \frac{\alpha}{\beta} \right); \left( 0,1 \right), \left( 1-\frac{1}{\alpha'}, \frac{2}{\alpha'} \right),\left( \frac{1}{2},1 \right)
    \end{array} \right]. \label{eq:38}
\end{eqnarray}
\normalsize
Thus, the energy eigenvalues are calculated using Eq. (\ref{eq:35}) as
\small
\begin{eqnarray}
    && E_\alpha = \frac{2 i^\alpha \eta A^2 t^{-\alpha-\frac{\alpha}{\beta}}} {\left( \alpha' \right)^2 \beta^2 \lambda^{\frac{2}{\alpha'}} \left( \frac{F^2}{\eta^2} \right)^{\frac{1}{2\beta}}} \nonumber \\
    && \times \int_0^{+\infty} H_{2,0:1,2;2,3}^{0,1:1,1;2,1}\left[ \left. 
    \begin{array}{c}
    \frac{-4\left( \frac{F}{\left(-i\right)^\alpha \eta} \right)^{\frac{1}{\beta}} t^{\frac{\alpha}{\beta}}} {x^2} \\
    \frac{\lambda^{-\frac{2}{\alpha'}}} {x^2}
    \end{array} \right|
    \begin{array}{c}
    \left( 1;1,1 \right), \left( \frac{1}{2};1,1 \right): \left( \frac{1}{2\beta}, \frac{1}{\beta} \right); \left( {1 - \frac{1}{\alpha'},\frac{2}{\alpha'}} \right),\left( 0,2 \right) \\
    - :\left( \frac{1}{2\beta}, \frac{1}{\beta} \right), \left(\frac{\alpha}{2\beta}, \frac{\alpha}{\beta} \right); \left( 0,1 \right), \left( 1-\frac{1}{\alpha'}, \frac{2}{\alpha'} \right),\left( \frac{1}{2},1 \right)
    \end{array} \right] \nonumber \\
    && \times H_{2,0:1,2;2,3}^{0,1:1,1;2,1}\left[ \left. 
    \begin{array}{c}
    \frac{-4\left( \frac{F}{i^\alpha \eta } \right)^{\frac{1}{\beta}} t^{\frac{\alpha}{\beta}}} {x^2} \\
    \frac{\lambda^{-\frac{2}{\alpha'}}} {x^2}
    \end{array} \right| 
    \begin{array}{c}
    \left( 1;1,1 \right), \left( \frac{1}{2};1,1 \right): \left( \frac{1}{2\beta}, \frac{1}{\beta} \right); \left( {1 - \frac{1}{\alpha'},\frac{2}{\alpha'}} \right),\left( 0,2 \right) \\
    - :\left( \frac{1}{2\beta}, \frac{1}{\beta} \right), \left( \alpha+\frac{\alpha}{2\beta}, \frac{\alpha}{\beta} \right); \left( 0,1 \right), \left( 1-\frac{1}{\alpha'}, \frac{2}{\alpha'} \right),\left( \frac{1}{2},1 \right)
    \end{array} \right] dx. \label{eq:39}
\end{eqnarray}
\normalsize
Using the identities of Fox's H-function, Eq. (\ref{eq:39}) can be written as
\small
\begin{eqnarray}
    && E_\alpha = \frac{2 i^\alpha \eta A^2 t^{-\alpha-\frac{\alpha}{\beta}}} {\left( \alpha' \right)^2 \beta^2 \lambda^{\frac{2}{\alpha'}} \left( \frac{F^2}{\eta^2} \right)^{\frac{1}{2\beta}}} \nonumber \\
    && \times \int_0^{+\infty} H_{0,2:2,1;3,2}^{1,0:1,1;1,2}\left[ \left. 
    \begin{array}{c}
    \frac{-x^2} {4\left( \frac{F}{\left(-i\right)^\alpha \eta} \right)^{\frac{1}{\beta}} t^{\frac{\alpha}{\beta}}} \\
    \lambda^{\frac{2}{\alpha'}}{x^2}
    \end{array} \right|
    \begin{array}{c}
    - : \left( 1-\frac{1}{2\beta},\frac{1}{\beta} \right), \left( 1-\frac{\alpha}{2\beta};\frac{\alpha}{\beta} \right) ; \left( 1,1 \right),\left( \frac{1}{\alpha'},\frac{2}{\alpha'} \right), \left( \frac{1}{2},1 \right) \\
    \left( 0;1,1 \right),\left( \frac{1}{2};1,1 \right):\left( 1-\frac{1}{2\beta},\frac{1}{\beta} \right);\left( \frac{1}{\alpha'},\frac{2}{\alpha'} \right),\left( 1,2 \right)
    \end{array} \right] \nonumber \\
    && \times H_{0,2:2,1;3,2}^{1,0:1,1;1,2}\left[ \left. 
    \begin{array}{c}
    \frac{-x^2} {4\left( \frac{F}{i^\alpha \eta} \right)^{\frac{1}{\beta}} t^{\frac{\alpha}{\beta}}} \\
    \lambda^{\frac{2}{\alpha'}}{x^2}
    \end{array} \right|
    \begin{array}{c}
    - : \left( 1-\frac{1}{2\beta},\frac{1}{\beta} \right), \left( 1-\alpha-\frac{\alpha}{2\beta};\frac{\alpha}{\beta} \right) ; \left( 1,1 \right),\left( \frac{1}{\alpha'},\frac{2}{\alpha'} \right), \left( \frac{1}{2},1 \right) \\
    \left( 0;1,1 \right),\left( \frac{1}{2};1,1 \right):\left( 1-\frac{1}{2\beta},\frac{1}{\beta} \right);\left( \frac{1}{\alpha'},\frac{2}{\alpha'} \right),\left( 1,2 \right)
    \end{array} \right] dx. \label{eq:40}
\end{eqnarray}
\normalsize
To calculate the integral (\ref{eq:40}) we first express both Fox's H-functions of two variables in terms of Mellin Barnes type integral as Eq. (\ref{eq:24}). So we have \cite{29}
\begin{eqnarray}
    && E_\alpha = \frac{2 i^\alpha \eta A^2 t^{-\alpha-\frac{\alpha}{\beta}}} 
    {\left( \alpha' \right)^2 \beta^2 \lambda^{\frac{2}{\alpha'}} \left( {\frac{F}{\eta}} \right)^{\frac{1}{\beta}}} \left( \frac{-1}{4\pi^2} \right) \int_{\mathcal{L}_1} \int_{\mathcal{L}_2} \int_{\mathcal{L}_3} 
    \frac{ \Gamma \left( {\frac{1}{\alpha'} - \frac{2}{\alpha'}t'} \right) \Gamma \left( t' \right) \Gamma \left( 1-\frac{1}{\alpha'} + \frac{2}{\alpha'}t' \right)} {\Gamma \left( {\frac{1}{2}-t'} \right) \Gamma \left( 2t' \right)} \left( \lambda^{\frac{2}{\alpha'}} \right)^{t'} \nonumber \\ 
    && \times \frac{\Gamma \left( 1-\frac{1}{2\beta} - \frac{1}{\beta} \xi' \right) \Gamma \left( \frac{1}{2\beta} + \frac{1}{\beta} \xi' \right)} {\Gamma \left( 1-\alpha-\frac{\alpha }{2\beta } - \frac{\alpha}{\beta}\xi' \right)} \left( \frac{-1}{4 \left( \frac{F}{{\left( i \right)}^\alpha \eta } \right)^{\frac{1}{\beta }} t^{\frac{\alpha}{\beta}}} \right)^{\xi'} \frac{\Gamma \left( \frac{1}{\alpha'} - \frac{2}{\alpha'} \eta' \right) \Gamma \left( \eta' \right) \Gamma \left( 1-\frac{1}{\alpha'} + \frac{2}{\alpha'} \eta' \right)} {\Gamma \left( {\frac{1}{2}-\eta'} \right) \Gamma \left( 2\eta' \right)} \nonumber \\
    && \times \left( \lambda^{\frac{2}{\alpha'}} \right)^{\eta'} \frac{\Gamma \left( -\xi'-\eta' \right)} {\Gamma \left( \frac{1}{2} + \xi' + \eta' \right)} 
    \left\{ \left[ \int_0^\infty \left( \frac{-1}{4\pi^2} \right) \int_{{\mathcal{L}_4}} \frac{\Gamma \left( -s-t' \right)} {\Gamma \left( \frac{1}{2}+s+t' \right)} \frac{\Gamma \left( 1-\frac{1}{2\beta} - \frac{1}{\beta} s \right) \Gamma \left( \frac{1}{2\beta} + \frac{1}{\beta} s \right)} {\Gamma \left( 1-\frac{\alpha}{2\beta} - \frac{\alpha}{\beta} s \right)} \right. \right. \nonumber \\ 
    && \times \left. \left. \left( \frac{-x^2} {4\left( \frac{F}{\left( -i \right)^\alpha \eta } \right)^{\frac{1}{\beta}} t^{\frac{\alpha}{\beta}}} \right)^s ds \right] x^{2\eta' + 2\xi' + 2t'} dx \right\} d\xi' dt' d\eta'. \label{eq:41}
\end{eqnarray}
Now if we use the Mellin transform theorem as if $f \left( x \right) = \frac{1}{2\pi i}\int_L x^{-s} g\left(s \right) ds$, then $g \left( s \right) = \int_0^\infty x^{s-1} f \left( x \right) dx$ and substitute $s \to -\eta' - t' - \xi' - \frac{1}{2}$ in Eq. (\ref{eq:41}), the following result is obtained
\footnotesize
\begin{eqnarray}
    && E_\alpha = \frac{2 i^\alpha \eta A^2 t^{-\alpha-\frac{\alpha}{\beta}}} 
    {\left( \alpha' \right)^2 \beta^2 \lambda^{\frac{2}{\alpha'}} \left( {\frac{F}{\eta}} \right)^{\frac{1}{\beta}}} \left( -4 \left( \frac{F}{\left( -i \right)^\alpha \eta } \right)^{\frac{1}{\beta }} t^{\frac{\alpha}{\beta}} \right)^{\frac{1}{2}} \nonumber \\
    && \times \int_{\mathcal{L}_1} \int_{\mathcal{L}_2} \int_{\mathcal{L}_3} 
    \frac{ \Gamma \left( {\frac{1}{\alpha'} - \frac{2}{\alpha'}t'} \right) \Gamma \left( t' \right) \Gamma \left( 1-\frac{1}{\alpha'} + \frac{2}{\alpha'}t' \right)} {\Gamma \left( {\frac{1}{2}-t'} \right) \Gamma \left( 2t' \right)} \frac{\Gamma \left( 1-\frac{1}{2\beta} - \frac{1}{\beta} \xi' \right) \Gamma \left( \frac{1}{2\beta} + \frac{1}{\beta} \xi' \right)} {\Gamma \left( 1-\alpha-\frac{\alpha }{2\beta } - \frac{\alpha}{\beta}\xi' \right)} \nonumber \\ 
    && \times \frac{\Gamma \left( 1 + \frac{\xi'}{\beta } + \frac{t'}{\beta} + \frac{\eta'}{\beta} \right) \Gamma \left( -\frac{{\xi'}}{\beta} - \frac{t'}{\beta} - \frac{\eta'}{\beta} \right)} {\Gamma \left( 1 + \frac{\alpha \xi'}{\beta} + \frac{\alpha t'}{\beta} + \frac{\eta \alpha'}{\beta} \right)}
    \left( -4 \lambda^{\frac{2}{\alpha'}} \left( \frac{F}{\left( -i \right)^\alpha \eta } \right)^{\frac{1}{\beta}} t^{\frac{\alpha}{\beta}} \right)^{t'} \nonumber \\
    && \times \left( -4 \lambda^{\frac{2}{\alpha'}} \left( \frac{F}{\left( -i \right)^\alpha \eta } \right)^{\frac{1}{\beta}} t^{\frac{\alpha}{\beta}} \right)^{\eta'}
    \left( -1 \right)^{\frac{\alpha \xi'}{\beta}} dt' d\xi' d\eta' \nonumber \\
    && = \frac{2 i^\alpha \eta A^2 t^{-\alpha-\frac{\alpha}{\beta}}} 
    {\left( \alpha' \right)^2 \beta^2 \lambda^{\frac{2}{\alpha'}} \left( {\frac{F}{\eta}} \right)^{\frac{1}{\beta}}} \left( -4 \left( \frac{F}{\left( -i \right)^\alpha \eta } \right)^{\frac{1}{\beta }} t^{\frac{\alpha}{\beta}} \right)^{\frac{1}{2}} 
    H_{1,2:2,1;3,2;3,2}^{1,1:1,1;1,2;1,2} \left[ \left.
    \begin{array}{c}
    \left( -1 \right)^{\frac{\alpha}{\beta}} \\
    -4\lambda^{\frac{2}{\alpha'}} \left( \frac{F}{\left( -i \right)^\alpha \eta } \right)^{\frac{1}{\beta}} t^{\frac{\alpha}{\beta}} \\
    -4 \lambda^{\frac{2}{\alpha'}} \left( \frac{F}{\left( -i \right)^\alpha \eta } \right)^{\frac{1}{\beta}} t^{\frac{\alpha}{\beta}} 
    \end{array} \right| \right. \nonumber \\
    && \left. \begin{array}{c}
    \left( 0;\frac{1}{\beta},\frac{1}{\beta},\frac{1}{\beta} \right):\left( 1-\frac{1}{2\beta},\frac{1}{\beta} \right),\left( 1 - \alpha  - \frac{\alpha }{{2\beta }},\frac{\alpha }{\beta } \right);\left( 1,1 \right),\left( \frac{1}{\alpha'},\frac{2}{\alpha'} \right),\left( \frac{1}{2},1 \right);\left( 1,1 \right),\left( \frac{1}{\alpha'},\frac{2}{\alpha'} \right),\left( \frac{1}{2},1 \right) \\
    \left( 0;\frac{1}{\beta},\frac{1}{\beta},\frac{1}{\beta} \right),\left( 0;\frac{\alpha}{\beta},\frac{\alpha}{\beta},\frac{\alpha}{\beta} \right):\left( 1-\frac{1}{2\beta},\frac{1}{\beta} \right);\left( \frac{1}{\alpha'},\frac{2}{\alpha'} \right),\left( 1,2 \right);\left( \frac{1}{\alpha'},\frac{2}{\alpha'} \right),\left( 1,2 \right)
    \end{array} \right], \label{eq:42}
\end{eqnarray}
\normalsize
where
$H_{p,q:\left[ {p_r},{q_r} \right]}^{m,n:\left[ {m_r},{n_r} \right]} \left[ \left. \begin{array}{c} x_1\\ \vdots \\ x_r \end{array} \right| \begin{array}{c}
\left( a_j;\alpha_j^1, \ldots ,\alpha_j^r \right)_{1,p}:\left\{ \left( {c_j^1, \gamma_j^1,} \right)_{1,p_r} \right\} \\
\left( b_j;\beta_j^1, \ldots ,\beta_j^r \right)_{1,q}:\left\{ \left( d_j^1,\delta_j^1, \right)_{1,q_r} \right\}
\end{array} \right]$
denotes Fox's H-function of several variables as follows \cite{37}
\begin{eqnarray}
    H_{p,q:\left[ {p_r},{q_r} \right]}^{m,n:\left[ {m_r},{n_r} \right]} \left[ \left. \begin{array}{c} x_1\\ \vdots \\ x_r \end{array} \right| \begin{array}{c}
    \left( a_j;\alpha_j^1, \ldots ,\alpha_j^r \right)_{1,p}:\left\{ \left( {c_j^1, \gamma_j^1,} \right)_{1,p_r} \right\} \\
    \left( b_j;\beta_j^1, \ldots ,\beta_j^r \right)_{1,q}:\left\{ \left( d_j^1,\delta_j^1, \right)_{1,q_r} \right\}
    \end{array} \right] \nonumber \\
    = \frac{1}{\left( 2\pi i \right)^r} \int_{\mathcal{L}_1} \ldots \int_{\mathcal{L}_r} \phi \left( \xi_1,\ldots,\xi _r \right) \prod_{i = 1}^r \left( \theta _i \left( \xi _i \right)x_i^{\xi_i} d\xi_i \right), \label{eq:43}
\end{eqnarray}
where
\begin{eqnarray}
    \phi \left( \xi_1, \ldots ,\xi_r \right) && = \frac{ \prod_{j = 1}^n \Gamma \left( 1 - a_j + \sum_{i = 1}^r \alpha_j^{\left( i \right)} \xi_i \right) \prod_{j = 1}^n \Gamma \left( b_j - \sum_{i = 1}^r \beta_j^{\left( i \right)} \xi_i \right)} { \prod_{j = n + 1}^p \Gamma \left( a_j - \sum_{i = 1}^r \alpha_j^{\left( i \right)} \xi_i \right) \prod_{j = 1}^n \Gamma \left( 1 - b_j + \sum_{i = 1}^r \beta_j^{\left( i \right)} \xi_i \right)}, \label{eq:44} \\
    \theta_i \left( \xi_i \right) && = \frac{ \prod_{j = 1}^{{n_i}} \Gamma \left( 1 - c_j^{\left( i \right)} + \gamma_j^{\left( i \right)} \xi_i \right) \prod_{j = 1}^{{m_i}} \Gamma \left( d_j^{\left( i \right)} - \delta_j^{\left( i \right)} \xi_i \right)} { \prod_{j = n_i + 1}^{p_i} \Gamma \left( c_j^{\left( i \right)} - \gamma_j^{\left( i \right)} \xi_i \right) \prod_{j = m_i + 1}^{q_i} \Gamma \left( {1 - d_j^{\left( i \right)} + \delta_j^{\left( i \right)} \xi_i} \right)}. \label{eq:45}
\end{eqnarray}
The integral Eq. (\ref{eq:42}) is valid under the following conditions
\begin{eqnarray*}
    \beta>-\frac{1}{4}, \frac{1}{2}<\alpha'<\frac{4}{3}, 2-\alpha>\beta.
\end{eqnarray*}
Thus, it is seen that the energy eigenvalue obtained in (\ref{eq:42}) depends on time, in other words the energy eigenvalue of Caputo space-time fractional Schr{\" o}dinger equation for the delta potential is time dependent. Now, we investigate the energy eigenvalue (\ref{eq:42}) for the special cases $\alpha \to 1$, $\beta \to 1$ and $\alpha' \to 1$ and then check the approximate behavior of energy at small and large times. Hence according to the above formula, we obtain
\begin{eqnarray}
    {E_\alpha}_{\left| \begin{array}{c}
    \alpha \to 1 \\
    \beta \to 1 \\
    \alpha' \to 1
    \end{array} \right.} 
    && = \frac{-2i \eta^2 A^2 t^{-2}} {\lambda^2 F \pi }
    \left( \frac{4Ft}{i\eta} \right)^{\frac{1}{2}} \int_{\mathcal{L}_1}  \int_{\mathcal{L}_2} \int_{\mathcal{L}_3} \Gamma \left( -\xi' - \eta' - t' \right) \Gamma \left( \xi' + \frac{3}{2} \right) \nonumber \\ 
    && \times \Gamma \left( t' \right) \Gamma \left( 1 - t' \right) \Gamma \left( \eta' \right) \Gamma \left( 1 - \eta' \right) \left( -1 \right)^{\xi'} \left( \frac{\lambda^2 Ft}{i \eta} \right)^{t'} \left( \frac{\lambda ^2 Ft}{i \eta } \right)^{\eta'} d\xi 'dt' d\eta' 
    = \ldots \nonumber \\ 
    && = \frac{ -4\left( \eta \right)^{\frac{3}{2}} A^2 (t)^{-\frac{3}{2}} (i)^{\frac{1}{2}}} {\lambda^2 \pi \left( F \right)^{\frac{1}{2}}} \sum_{n = 0}^\infty \frac{i}{n!} H_{1,0:1,1;1,1}^{0,1:1,1;1,1}\left[ \left. \begin{array}{c}
    \frac{i\eta}{\lambda^2 Ft} \\
    \frac{i\eta}{\lambda^2 Ft}
    \end{array} \right| \begin{array}{c}
    \left( - \frac{1}{2} - n \right) : \left( 0,1 \right) ; \left( 0,1 \right) \\
    - :\left( 0,1 \right) ; \left( 0,1 \right)
    \end{array} \right] \nonumber \\ 
    && = \frac{ -4 \left( \eta \right)^{\frac{3}{2}} A^2 (t)^{-\frac{3}{2}} (i)^{\frac{1}{2}}} {\lambda^2\pi \left( F \right)^{\frac{1}{2}}} \sum_{n = 0}^\infty \frac{i}{n!}H_{1,2}^{2,1}\left[ \frac{\lambda^2 Ft} {i\eta}\left| \begin{array}{c}
    \left( 1,1 \right) \\
    \left( \frac{3}{2} + n,1 \right),\left( 2,1 \right)
    \end{array} \right. \right]. \label{eq:46}
\end{eqnarray}
In the above relation, we have used the identities of Fox's H-function of two variables given in Refs. \cite{32,33,34}. We also need to use the expansion of the Fox's H-function at small and large times. For $t \to 0$, using relations of Fox's H-function \cite{34}, Eq. (\ref{eq:46}) is obtained as 
$ E_\alpha = \frac{-4 \left( \eta \right)^{\frac{3}{2}} A^2 (t)^{-\frac{3}{2}} (i)^{\frac{1}{2}}} {\lambda^2 \pi \left( F \right)^{\frac{1}{2}}} \times \frac{i\pi}{4} \times \left( \frac{\lambda^2 Ft}{i\eta} \right)^{\frac{3}{2}} = -A^2 F\lambda  = - \frac{m v_0^2}{2\hbar^2}$, and it is seen that the energy value is in accordance with the eigenvalue for the ordinary Schr{\" o}dinger equation. For $t \to \infty$ and by using the Fox's H-function properties \cite{34}, we get
$E_\alpha = \frac{-4\left( \eta \right)^{\frac{3}{2}} A^2 (t)^{ -\frac{3}{2}} (i)^{\frac{1}{2}}} {\lambda^2 \pi \left( F \right)^{\frac{1}{2}}}$. This phrase tends to zero in $t \to \infty$, that is, the energy value goes to zero.

\section{\label{section5}Conclusions}
We have investigated Caputo space-time fractional Schr{\" o}dinger equation with the delta potential for the initial and boundary imposed condition. In this condition, we have obtained the wave function and the energy eigenvalues, and it is seen that the results are consistent with the ordinary Schr{\" o}dinger equation for $\alpha \to 1$, $\beta \to 1$ and $\alpha' \to 1$ cases.

\end{document}